\documentclass[12pt,psfig]{amsart}
\usepackage{amsmath,amsthm,epsfig,amssymb}
\def\ra#1{\hbox to #1pc{\rightarrowfill}}
\def\fract#1#2{\raise4pt\hbox{$ #1 \atop #2 $}}

\newtheorem{proposition}{Proposition}
\newtheorem{theorem}{Theorem}
\newtheorem{lemma}{Lemma}
\newtheorem{corollary}{Corollary}
\newtheorem{remark}{Remark}
\newtheorem{example}{Example}

\newtheorem{definition}{Definition}

\newenvironment{pf}{{\bf Proof:\
}}{\raisebox{-4pt}{\rule{6pt}{6pt}}\vspace{4pt plus 8pt minus 1pt}\vskip
8pt}

\newenvironment{prooft}{{\leftmargin=20pt\rightmargin=0pt} {\bf Proof of
Theorem:\ }}{\vspace{4pt plus 2pt minus 1pt}}

\newcommand{\bth}{\begin{theorem}}
\newcommand{\bpr}{\begin{proposition}}
\newcommand{\epr}{\end{proposition}}
\newcommand{\bco}{\begin{corollary}}
\newcommand{\eco}{\end{corollary}}
\newcommand{\ble}{\begin{lemma}}
\newcommand{\ele}{\end{lemma}}
\newcommand{\bre}{\begin{remark}\rm}
\newcommand{\ere}{\end{remark}}
\newcommand{\bex}{\begin{example}\rm}
\newcommand{\eex}{\end{example}}
\newcommand{\bres}{{\bf Remarks}\begin{enumerate}}
\newcommand{\eres}{\end{enumerate}}
\newcommand{\bexs}{{\bf Examples}\begin{enumerate}\begin{enumerate}}
\newcommand{\eexs}{\end{enumerate}\end{enumerate}}
\newcommand{\bde}{\begin{definition}\rm}
\newcommand{\ede}{\end{definition}}
\newcommand{\bpf}{\begin{pf}}
\newcommand{\epf}{\end{pf}}
\newcommand{\bproof}{\begin{pf}}
\newcommand{\eproof}{\end{pf}}
\newcommand{\bpft}{\begin{prooft}\begin{itemize}}
\newcommand{\epft}{\end{itemize}\end{prooft}}
\begin{document}
\date{\today}

\title{\bf GENUS ZERO  ACTIONS ON RIEMANN SURFACES}
\author{\bf SADOK KALLEL AND DENIS SJERVE}
\thanks{Research supported by PIms 
NSERC grant A7218}
\maketitle
\noindent\section{\bf INTRODUCTION}\label{intro}

In this paper we solve the following problem:

\noindent {\bf PROBLEM.}  Which finite groups $G$ admit an action on
some compact connected Riemann surface $M$ so that if $H$ is any
non-trivial subgroup of $G$ then the orbit surface $M/H$ has genus
zero, that is
\begin{eqnarray}\displaystyle\label{genuscond1}
1\ne H\subseteq G\Longrightarrow M/H=\mathbb{P}^1(\mathbb{C}). 
\end{eqnarray}
Any finite group $G$ will admit infinitely many actions $G\times M\to
M$ so that $M/G=\mathbb P^1(\mathbb C).$ It will turn out that very
few groups $G$ satisfy (\ref{genuscond1}).

The action $G\times M\to M$ is supposed to be analytic and effective.
Thus $G$ is a subgroup of $\mathrm{Aut}(M)$, the group of all analytic
automorphisms of $M.$ For any action $G\times M\to M$, if
$\displaystyle M/H_1=\mathbb{P}^1(\mathbb{C})$ for some $\displaystyle
H_1\subseteq G$, then automatically $\displaystyle
M/H_2=\mathbb{P}^1(\mathbb{C})$ if $\displaystyle H_1\subseteq
H_2\subseteq G$.  This means that we can reduce the problem to the
consideration of cyclic subgroups $\mathbb Z_p\subseteq G$, where $p$
is a prime dividing the order of $G$.

Our problem thus becomes:

\noindent {\bf PROBLEM.}  Which finite groups $G$ admit an action on
some compact connected Riemann surface $M$ so that if $H$ is any
cyclic subgroup of prime order $p$ then the orbit surface $M/H$ has
genus zero, that is
\begin{eqnarray}\displaystyle\label{genuscond2}
M/\mathbb Z_p=\mathbb{P}^1(\mathbb{C}),~\mathrm{for~all}~ Z_p\subseteq G,~
\mathrm{where}~p~\mathrm{is~a~prime~dividing}~|G|.
\end{eqnarray}
All Riemann surfaces considered in this paper, with just a few
exceptions, will be compact and connected.  The exceptions are the
complex plane and the upper half plane. The groups being studied in
this paper will always be finite.

In our solution of this problem not only do we describe the groups $G$
admitting such actions, but we also determine all possible actions for
each group.  This amounts to describing all admissible epimorphisms
$\theta:\Gamma\to G,$ where $\Gamma$ is a Fuchsian group of signature
$(0~|~n_1,\ldots,n_r).$ See Section (\ref{prelim}).

In Section (\ref{g=0 or 1}) we consider the low genera cases $g=0,1,$
where $g$ is the genus of $M.$ Other than these exceptions we will
usually assume that $g>1.$

\bde 
A group $G$ is said to have genus zero if there exists a Riemann
surface $M$ and an action of $G$ on $M$ satisfying (\ref{genuscond1}),
or equivalently (\ref{genuscond2}) . We also say that the action has
genus zero.  
\ede

To explain the significance of the above problem, and to give it some
motivation, we recall the definition of a fixed point free linear
action on a $\mathbb C$ vector space $V$.  
\bde 
A linear action $G\times V\to V$ is said to be fixed point free if
$$
S\in G,~S\ne 1\Longrightarrow~S(v)\ne v~\mathrm{for~all}~v\in V,~v\ne 0.
$$
\ede
Now let $V$ be the vector space of holomorphic differentials on
$M$. $V$ is a $\mathbb C$ vector space of dimension $g$, where $g$ is
the genus of $M$. The action of $G$ on $M$ induces a linear action on
the vector space $V.$ If $S\in G$ is an element of order $p$ then the
induced linear transformation $S^*:V\to V$ also has order $p$,
assuming $g>1.$ It follows that the eigenvalues of $S^*$ are $p^{th}$
roots of unity.  In particular $+1$ might be an eigenvalue of
$T^*$. Part of the proof of the Eichler trace formula, see Farkas and
Kra \cite{Farkas}, states that the dimension of the $+1$ eigenspace of
$S^*$ is the genus of $M/H$, where $H\cong \mathbb Z_p$ is the cyclic
group of order $p$ generated by $S.$

Therefore the induced action $G\times V\to V$ is fixed point free if,
and only if, $G\times M\to M$ is a genus-zero action. We can impose a
metric on $V$ so that the action becomes unitary. Thus the action
$G\times M\to M$ has genus zero if, and only if, $S^{2g-1}/G$ is an
elliptic space form. This equivalence assumes that $g>1.$ Here
$S^{2g-1}$ denotes the unit sphere in $V$.

The converse problem of when a fixed point free linear action $G\times
V\to V$ arises from a genus-zero action on some Riemann surface is not
considered in this paper.

Groups $G$ which admit fixed point free linear actions have been
classified, see Wolf \cite{Wolf}.  In particular they must satisfy a
strong condition on their Sylow $p$-subgroups. First recall that the
generalized quaternion group $Q(2^n)$ is defined as follows: 
\bde 
The generalized quaternion group is the group with the presentation:
\begin{equation}\displaystyle\label{presQ}
Q(2^n)=\left<A,B~|~A^{2^{n-1}}=1,~B^2=A^{2^{n-2}},~BAB^{-1}=A^{-1}\right>.
\end{equation}
\ede 
We will always assume $n\ge 3,$ since otherwise $Q(2^n)$ is cyclic.  
\bde \label{sylowcond} 
We say that a group $G$ satisfies the { \bf Sylow conditions} if the
following two conditions hold.
\begin{enumerate}
\item For an odd prime $p$ the Sylow $p$-subgroups are cyclic. 
\item The Sylow $2$-subgroups are either cyclic or generalized quaternion.
\end{enumerate}
\ede
The Sylow conditions are equivalent to the following:
\bde\label{p^2cond}
We say that a group $G$ satisfies the {\bf ${\bf p^2}$ conditions}
if every subgroup of order $p^2$ is cyclic,
where $p$ is any prime.
\ede 

Every group $G$ admitting a fixed point free linear action satisfies
these conditions.  In fact these groups must satisfy the even stronger
$pq$ conditions.
\bde\label{pqcond}
A group $G$ satisfies the {\bf ${\bf pq}$ conditions} if  every subgroup
of order $pq$ is cyclic, where 
$p$ and $q$ are arbitrary primes.
\ede
Let $D_{2n}$ denote the dihedral group of order $2n$. If $n$ is odd
then $D_{2n}$ satisfies the $p^2$ conditions but not the $pq$
conditions. In fact any group $G$ of even order which admits a fixed
point free linear action must have exactly one element of order $2,$
and this element generates the center of $G.$

To describe our results we need another definition and some notation.
\bde\label{ZMgrps} We let
$G_{m,n}(r)$ denote the group presented as follows:
\begin{eqnarray}\displaystyle\label{ZMpres}
\notag\mathrm{generators:}& &A,B;\\
\mathrm{~relations:}& & A^m=1,~B^n=1,~BAB^{-1}=A^r;\\
\notag
\mathrm{conditions:}& &GCD((r-1)n,m)=1~\mathrm{and}~r^n\equiv 1~(mod~m).
\end{eqnarray} 
\ede
These groups are precisely the groups having all Sylow subgroups
cyclic, see Burnside \cite{Burnside}.  To avoid the trivial cases
where the group is cyclic we will usually assume that $m>1,n>1.$ Note
that the conditions imply $r\not\equiv 1~(mod~m).$

If $d$ denotes the order of $r$ modulo $m$ then Zassenhaus
\cite{Zassenhaus} proved that $G_{m,n}(r)$ satisfies the $pq$
conditions if, and only if, every prime divisor of $d$ also divides
$\displaystyle\frac{n}{d}.$ He also proved that $G_{m,n}(r)$ admits a
fixed point free linear representation if, and only if, the $pq$
conditions hold.  The groups $G_{m,n}(r)$ are known as Zassenhaus
metacyclic groups (abbreviated to ZM groups).

Let $I^*$ denote the binary icosahedral group. It has order $120$ and
admits a fixed point free linear representation. $I^*$ is
non-solvable, and if $G$ is a non-solvable group admitting a linear
fixed point free representation then $G$ contains $I^*$ as a subgroup.

The main results of this paper are contained in the following
theorems. In particular the first theorem implies that $I^*$ does not
have genus zero, and therefore neither does any non-solvable group
admitting a fixed point free linear representation.
\begin{theorem}
The  groups having genus zero are the cyclic groups, the generalized
quaternion groups $Q(2^n),$ the polyhedral
groups and the ZM groups $G_{p,4}(-1),$ where $p$ is an odd 
prime.
\end{theorem}\label{mainth}

The cyclic and polyhedral groups have genus zero because they can act
on $\mathbb P^1(\mathbb C).$ Some cyclic groups admit genus zero
actions on surfaces of higher genus, but most do not. See Theorems
(\ref{cyclicgenusactions}), (\ref{pqactions}), (\ref{p^2qactions}),
and Corollary (\ref{cyclicactions}).
 
This theorem gives a solution to Problem (\ref{genuscond1}) but does
not describe the actions involved. To do this we need some more
notation. Let $\Gamma(0~|~n_1,\ldots,n_r)$ denote the abstract group
presented by
\begin{eqnarray}\label{gppres}
\left< X_1,X_2,\ldots,X_r \bigm |
 X_{1}^{n_1}=X_{2}^{n_2}=\cdots=X_{r}^{n_r}=
        X_{1}X_{2}\cdots X_{r}=1\right>.
\end{eqnarray}
For a more detailed explanation of the notation see
Section (\ref{prelim}).
  
If $G\times M\to M$ is an action satisfying
$M/G=\mathbb{P}^1(\mathbb{C}),$ and the genus of $M$ is $g,$ then
there exists a short exact sequence
\begin{eqnarray}\label{ses1}
1\to \Pi\to \Gamma\stackrel{\theta}{\longrightarrow}G\to 1,
\end{eqnarray}
where $\Pi\cong\pi_1(M)$ and $\Gamma=\Gamma(0~|~n_1,\ldots,n_r)$ 
for some choice of $n_j.$ We say that the signature of the action is
$(0~|~n_1,\ldots,n_r).$ The signature of the action, together with the
epimorphism $\theta:\Gamma\to G$ and the particular realization of $\Gamma$
as a Fuchsian group, completely determines the action.

\begin{theorem}$\displaystyle\label{cyclicgenusactions}$
All  genus-zero actions of the cyclic group $\mathbb Z_{p^e},$ where
$p\ge 2$ is any prime, have signature 
$(0|\underbrace{p,\ldots,p}_{r},p^e,p^e),$ where $r$ is arbitrary. 
The genus is
$g=\frac{1}{2}r\left(p^e-p^{e-1}\right).$
\end{theorem}
If $e=1$ this means that the signature of the action is
$(0|\underbrace{p,\ldots,p}_{r})$, where $r\ge 2,$ and the genus is
$g=\frac{1}{2}(r-2)(p-1).$

\begin{theorem}\label{quatgenusactions}
$\displaystyle$
All genus-zero actions of the generalized quaternion group 
$Q(2^n)$ have signature 
$(0|\underbrace{2,\ldots,2}_{r},4,4,2^{n-1}),$ where $r$ is odd. The genus
is
$g=2^{n-2}(r+1).$
\end{theorem}

Theorem 2 shows
there are infinitely many signatures for genus-zero actions by
cyclic $p-$groups. This is not true for other cyclic groups.
\begin{theorem}\label{pqactions}
Suppose $p,q$ are distinct primes. Then the  genus-zero actions of 
$\mathbb Z_{pq}$ have signature and corresponding genus given by
\begin{enumerate}
\item sig$(\Gamma)=(0~|~pq,pq)$, in which case $g=0.$
\item sig$(\Gamma)=(0~|~p,q,pq)$, in which case 
$
\displaystyle
g=\frac{1}{2}(p-1)(q-1).
$
\item[3.] sig$(\Gamma)=(0~|~p,p,q,q)$, in which case $g=(p-1)(q-1).$
\end{enumerate}
\end{theorem}

See Section (\ref{sylow}) for the proofs of Theorems
(\ref{cyclicgenusactions}), (\ref{quatgenusactions}) and Section
(\ref{cyclicgrps}) for the proof of Theorem (\ref{pqactions}).  These
sections also give details of the epimorphisms $\theta:\Gamma\to G$
classifying the actions.

\begin{theorem}\label{ZMactions}
All genus-zero actions of the ZM group $G_{p,4}(-1),$ where $p$
is an odd prime, 
have  signature $(0~|~4,4,p).$  The corresponding genus is $g=p-1.$
\end{theorem}
This theorem is proved in Section (\ref{zassgrps}). In Section 
(\ref{binaryicos}) we complete the proof of Theorem (\ref{mainth}).

\section{PRELIMINARIES}\label{prelim}

In this section we collect some preliminary material and review some
of the material on Riemann surfaces that we need later.

Let $\mathbb U$ denote a simply connected Riemann surface, that is
$\mathbb P^1(\mathbb C),$ $\mathbb C$ or the upper half plane $\mathbb
H$.  A Fuchsian group $\Gamma$ is any finitely generated discrete
subgroup of $PSL_2(\mathbb R)$, the group of analytic automorphisms of
$\mathbb H$. By abuse of terminology we shall also call a finitely
generated discrete subgroup of Aut$(\mathbb U)$ Fuchsian.

To every Fuchsian group $\Gamma$ we associate a signature
$(h~|~n_1,n_2,\ldots,n_r)$, where $h$ is the genus of the orbit
surface $\mathbb U/\Gamma$ and $n_1,\ldots,n_r$ are the orders of the
distinct conjugacy classes of maximal cyclic subgroups of $\Gamma.$
The $n_j$ are called the periods.  In general it is possible that some
of the $n_j$ are infinite, but the Fuchsian groups $\Gamma$ considered
in this paper will all have the property that $\mathbb U/\Gamma$ is
compact, and so the periods $n_j$ will be finite. We use the notation
$\mathrm{sig}(\Gamma)$ for the signature of $\Gamma.$

The notation $\Gamma(h~|~n_1,\ldots,n_r)$ denotes any Fuchsian group
of signature $(h~|~n_1,n_2,\ldots,n_r)$. If $r=0$ this group is
torsion free and we use the notation $\Gamma(h~|-).$ In fact the
Fuchsian groups with $r=0$ are just the fundamental groups of Riemann
surfaces.

Of particular interest to us will be those Fuchsian groups
$\Gamma=\Gamma(0~|~n_1,\ldots,n_r)$ since they play a seminal role in
actions $G\times M\to M$ satisfying $M/G=\mathbb P^1(\mathbb C).$ As
an abstract group $\Gamma$ has the presentation (\ref{gppres}).  The
geometry of $\Gamma$ is spherical, euclidean or hyperbolic according
as $\sum_{j=1}^{r}\frac{1}{n_j}>1,~=1$ or $<1.$ In the hyperbolic
case, for any realization of $\Gamma$ as a Fuchsian group, the $X_j$
are elliptic and so are rotations about vertices $V_j\in\mathbb H.$ Up
to conjugation by elements of Aut$(\mathbb U)$ the space of all such
realizations is a cell of dimension $2r-6.$

Let $G$ be a group acting on a Riemann surface $M$ of genus $g$ and
let $\mathbb U$ be the universal covering space of $M.$ Then there
exists a short exact sequence of groups
\begin{eqnarray}
1\to \Pi\to \Gamma\stackrel{\theta}{\longrightarrow}
G\to 1,
\end{eqnarray}
where 
\begin{enumerate}
\item $\Gamma$ is  a Fuchsian group  with signature 
$(h~|~n_1,n_2,\ldots,n_r).$
\item
$\Pi=\mathrm{Ker}(\theta)$ is a torsion free Fuchsian group with signature
$(g~|-)$.
\item
$M=\mathbb U/\Pi$ and the action of an element $S\in G$ on $M$ is given by 
$S[z]=[\gamma(z)]$, where the brackets
[ ] indicate the $\Pi$ equivalence class of points in $\mathbb U$ and
$\gamma\in\Gamma$ is any element such that $\theta(\gamma)=S.$ 
\item The orbit surface $M/G$ has genus $h$ and is naturally isomorphic to 
$\mathbb U/\Gamma.$
\end{enumerate}

The relationship between the genera $g,h$ is given by the
Riemann-Hurwitz formula:
\begin{equation}\label{RH1}\displaystyle
2g-2=|G|\left(2h-2+\sum_{j=1}^{r}\left(1-\frac{1}{n_j}\right)\right).
\end{equation}
\bde
The signature of the action of $G$ on $M$ is defined to be the 
signature of $\Gamma.$
\ede

Suppose $S\in G$ is an element of order $p$, where $p$ is a prime, and
$H$ is the cyclic subgroup generated by $S.$ Then
$\Gamma^{'}=\theta^{-1}(H)$ is a Fuchsian group and there exists a
short exact sequence
\begin{equation}\label{ses2}
1\to \Pi\to \Gamma^{'}\stackrel{\theta}{\longrightarrow}
 H\to 1.
\end{equation}
The signature of $\Gamma^{'}$ will have the form
$(k~|~\underbrace{p,\ldots,p}_{t}),$ where $t$ is the number of fixed
points of $S:M\to M.$ The Riemann-Hurwitz formula (\ref{RH1}) gives
\begin{eqnarray}\displaystyle\label{RH2}
2g-2=p\left(2k-2+t\left(1-\frac{1}{p}\right)\right).
\end{eqnarray}
A comparison between Formulas (\ref{RH1}) and (\ref{RH2}) then gives a
numerical restriction on the possible genera and the number of fixed
points.

If $G\times M\to M$ is an action satisfying $M/G=\mathbb P^1(\mathbb
C)$ then the signature of $\Gamma$ will be given by
$\mathrm{sig}(\Gamma)=(0~|~n_1,\ldots,n_r),$ and so $\Gamma$ will be
given abstractly by the presentation in (\ref{gppres}).  In this case
let
$$T_j=\theta(X_j)~\mathrm{and}~G_j=~\mathrm{the~subgroup~of}~G~
\mathrm{generated~by}~ T_j,1\le j\le r.$$ Then 
\begin{enumerate}

\item  $T_1,T_2,\ldots,T_r$ generate $G$ ( $\theta$ is an epimorphism).

\item 
$T_{1}^{n_1}=T_{2}^{n_2}=\cdots= T_{r}^{n_r}=T_1T_2\cdots T_r=1$ ($\theta$
must preserve the relations in $\Gamma$).

\item The order of $T_j$ is  $n_j$, $1\le j\le r$ (the kernel of
$\theta$ is torsion free).

\end{enumerate}

The converse is true. In other words if we are given
$T_1,T_2,\ldots,T_r$ satisfying these conditions then there is an
action $G\times M\to M$ satisfying $M/G=\mathbb P^1(\mathbb C)$ and
having signature $(0~|~n_1,\ldots,n_r).$ It follows that every group
$G$ admits infinitely many actions $G\times M\to M$ satisfying
$M/G=\mathbb P^1(\mathbb C).$

\bre\label{fixedpts1}
The key to the study of actions $G\times M\to M$ satisfying
$M/G=\mathbb P^1(\mathbb C)$ is the determination of the number of
fixed points of the action. Let $S\in G$, $S\ne 1,$ and choose any
$\gamma\in \Gamma$ such that $\theta(\gamma)=S.$ Then a point $[z]\in
M$ will be a fixed point of $S$ if, and only if, there exists $\mu\in
\Pi$ such that $\mu\gamma(z)=z.$ Now the elements of $\Gamma$ that
have fixed points are the conjugates of powers of the elliptic
generators.  That is we must have
$$\displaystyle
\mu\gamma
=\delta X_{j}^{k_j}\delta^{-1},~\mathrm{where}~1\le j\le r,~1\le k_j <n_j,~
\delta\in\Gamma,
$$
in which case $z=\delta(V_j).$  Thus the fixed points
of $S$ are those $\Pi$ equivalence classes $[\delta(V_j)]$, any 
$\delta\in \Gamma$, satisfying
\begin{eqnarray}\displaystyle\label{cond}
S=dT_{j}^{k_j}d^{-1},~\mathrm{for~some}~k_j,~1\le j\le r,~
\mathrm{where}~d=\theta(\delta).
\end{eqnarray}
\ere
Condition (\ref{cond}) is especially easy to use if
$S\in G,S\ne 1,$ is in the center.

\ble\label{centralelt}
Suppose $G$ is a group admitting an action $G\times M\to M$ satisfying
$M/G=\mathbb P^1(\mathbb C)$ and let $S\in G,S\ne 1,$ be a central
element.  Then the fixed points of $S$ are the $\Pi$ equivalence
classes $[\delta(V_j)],~ 1\le j\le r,$ any $\delta\in \Gamma$, such
that $S\in G_j$.  \ele

\bco\label{fixedpts2}
Suppose $G$ is a group admitting an action $G\times M\to M$ satisfying
$M/G=\mathbb P^1(\mathbb C)$ and let $S\in G$ be a central element of
order $p.$ Then the number of fixed points of $S$ is $\displaystyle
|G|{\sum_{j}}^{'}{\frac{1}{n_j}}$, where the prime indicates we sum
only over those $j$ such that $S\in G_j.$ \eco \bpf According to Lemma
(\ref{centralelt}) the fixed points of $S$ are those $\Pi$ equivalence
classes $[\delta(V_j)],~ 1\le j\le r$, any $\delta\in\Gamma,$ such
that $S\in G_j$.  If $S\in G_j\cap G_k$ then
\begin{eqnarray*}\displaystyle
[\delta(V_j)]=[\epsilon(V_k)]&\Longleftrightarrow&
\mu\delta(V_j)=\epsilon(V_k)~\mathrm{for~some}~\mu\in\Pi\\
&\Longleftrightarrow& j=k~\mathrm{and}~
\theta(\delta)\equiv\theta(\epsilon)(mod~G_j).
\end{eqnarray*}
Thus the fixed points of $S$, for a fixed $j$ with $S\in G_j$, are in 
one-to-one
correspondence with the cosets of $G_j$ in $G.$
\epf
In a similar fashion we can prove the next corollary.

\bco\label{fixedpts3}
Suppose $G$ is a group admitting an action $G\times M\to M$ satisfying
$M/G=\mathbb P^1(\mathbb C).$ Suppose $G$ has a unique subgroup of
order $p,$ where $p$ is a prime, and let $S$ be any element of order
$p$.  Then the number of fixed points of $S$ is $\displaystyle
|G|{\sum_{j}}^{'}{\frac{1}{n_j}},$ where the prime indicates we sum
over those $j$ so that $p|n_j.$ \eco

The cyclic groups $\mathbb Z_{p^{e}}$ have unique subgroups of orders
$p,p^2,\ldots,p^e$ and the generalized quaternion goup $Q(2^n)$ has a
unique subgroup of order $2,$ the subgroup generated by the central
element $B^2.$ In fact the only $p$-groups which contain a unique
subgroup of order $p$ are the cyclic groups and the generalized
quaternion groups.

Let $G$ denote either of these groups and suppose $S\in G$ is an
element of order $p.$ If $G$ acts on a Riemann surface $M$ so that
$M/G=\mathbb P^1(\mathbb C)$ then there is a short exact sequence as
in (\ref{ses1}), and moreover the signature of $\Gamma$ must have the
form
\begin{eqnarray}\label{cyclicsig}
\mathrm{sig}(\Gamma)&=&(0|\underbrace{p,\ldots,p}_{r_1},
\underbrace{p^2,\ldots,p^2}_{r_2},
\ldots,
\underbrace{p^e,\ldots,p^e}_{r_e})~\mathrm{if}~ G=\mathbb Z_p^{e}. \\
\label{quatsig}
\mathrm{sig}(\Gamma)&=&(0|\underbrace{2,\ldots,2}_{r_1},
\underbrace{4,\ldots,4}_{r_2},
\ldots,
\underbrace{2^{n-1},\ldots,2^{n-1}}_{r_{n-1}})~\mathrm{if}~G=Q(2^n).
\end{eqnarray}
Not all signatures are realizable since there are restrictions on the
$r_j$ that must be satisfied in order that $\theta:\Gamma\to G$ be
well defined, onto and have torsion free kernel. See Section
(\ref{sylow}).  

\bco\label{fixedpts4} Let $G$ denote either $\mathbb Z_{p^e}$ or
$Q(2^n)$, and $S$ an element of order $p$. Suppose $G$ acts on a
Riemann surface $M$ so that $M/G=\mathbb P^1(\mathbb C).$ Then the
number of fixed points of $S$ is
$$\displaystyle p^e\sum_{j=1}^{e}\frac{r_j}{p^j}=
\sum_{j=1}^{e}r_jp^{e-j}~\mathrm{if}~G=\mathbb Z_{p^e}~\mathrm{and}~
\displaystyle 2^n\sum_{j=1}^{n-1}\frac{r_j}{2^j}
=\sum_{j=1}^{n-1}r_j2^{n-j}~\mathrm{if}~G=Q(2^n).$$

\eco
\begin{proof}
This is an immediate consequence of the last corollary.
\end{proof}

\section{ACTIONS ON A SURFACE OF GENUS 0 OR 1}\label{g=0 or 1}

In this section we determine all genus-zero actions on either $\mathbb
P^1(\mathbb C)$ or a torus.  Such actions on $\mathbb P^1(\mathbb C)$
are very easy to describe since the possible groups are just the
finite subgroups of $PSL_2(\mathbb C).$ These are well known to be
either cyclic or polyhedral. In terms of our notation for Fuchsian
groups the finite subgroups of $PSL_2(\mathbb C)$ are:
\begin{enumerate}
\item $\Gamma(0~|~n,n),$ the cyclic group of order $n.$
\item $\Gamma(0~|~2,2,n),$ the dihedral group of order $2n.$
\item $\Gamma(0~|~2,3,3),$ the tetrahedral group of order $12.$
\item $\Gamma(0~|~2,3,4),$ the octahedral group of order $24.$
\item $\Gamma(0~|~2,3,5),$ the icosahedral group of order $60.$ 
\end{enumerate} 

This also gives the possible signatures of the actions.  Moreover, up
to conjugacy in $PSL_2(\mathbb C),$ each of these groups admits a
unique embedding and so a unique genus-zero action on $\mathbb
P^1(\mathbb C).$

Now suppose $G$ acts on a torus $M$ so that (\ref{genuscond1}) is
satisfied.  It is well known that the group of automorphisms Aut($M$)
contains $M$ as a normal subgroup, acting on itself by translations,
and that the quotient group is cyclic of order $2,4$ or $6.$ In other
words there is a short exact sequence
$$
1\to M\to \mbox{Aut}(M)\to K\to 1,~\mbox{where}~ 
K\cong \mathbb Z_2,\mathbb Z_4,~\mbox{or}~\mathbb Z_6.
$$ 
The elements of the subgroup $M$ act fixed point freely and therefore
finite subgroups of $M$ will not satisfy (\ref{genuscond1}).  In fact
the quotient by any finite subgroup will again be a torus.  It follows
that $G$ is a subgroup of $\mathbb Z_2,\mathbb Z_4$ or $\mathbb Z_6,$
and so $G\cong \mathbb Z_2,\mathbb Z_3,\mathbb Z_4$ or $\mathbb Z_6.$

The next theorem summarizes the situation.
\begin{theorem}\label{g=1}
The groups acting on a torus $M$ and satisfying (\ref{genuscond1}) 
are given as follows:
\begin{eqnarray*}
& &\mathbb Z_2~\mathrm{with~signature}~(0~|~2,2,2,2).\\
& &\mathbb Z_3~\mathrm{with~signature}~(0~|~3,3,3).\\
& &\mathbb Z_4~\mathrm{with~signature}~(0~|~2,4,4).\\ 
& &\mathbb Z_6~\mathrm{with~signature}~(0~|~2,3,6).
\end{eqnarray*}
\end{theorem}
In all cases the epimorphism $\theta:\Gamma\to G$ is unique up to
automorphisms of $G.$
\section{SYLOW GROUP ACTIONS}\label{sylow}

In this section we prove Theorems (\ref{cyclicgenusactions}) and
(\ref{quatgenusactions}), that is we classify all genus-zero actions
of the cyclic groups $\mathbb Z_{p^e}$, where $p\ge 2$ is any prime,
and the generalized quaternion groups $Q(2^n).$ For each of these
groups there are infinitely many possible signatures.

First we consider the case $G=\mathbb Z_{p^e}$. Let $T$ be a generator
of $G$ and suppose $G$ acts on a Riemann surface $M$ so that
$M/G=\mathbb P^{1}(\mathbb C)$. For the moment we do not assume that
the action has genus zero.  Then there exists a short exact sequence
as in (\ref{ses1}) and sig($\Gamma$) is given by (\ref{cyclicsig}).

We choose elliptic generators of $\Gamma$ as follows:
$$\displaystyle
X_{j,k}~\mathrm{corresponding~to~the~period}~p^j,~1\le j\le e,~1\le k\le r_j.
$$

Then the epimorphism $\theta :\Gamma\to\mathbb Z_{p^e}$ is given by:
\begin{equation}\displaystyle\label{eq1}
\theta(X_{j,k})=T^{a_{j,k}p^{e-j}},~\mathrm{where}~1\le j\le e,~
1\le k\le r_j~\mathrm{and}~a_{j,k}\not\equiv 0~(mod~p).
\end{equation}

The conditions on the $a_{j,k}$ guarantee that $\Pi$= Ker$(\theta)$ is
torsion free.  To ensure that $\theta$ is onto and well defined we
also need the conditions:
\begin{equation}\displaystyle\label{eq2}
r_e\ge 2~\mathrm{and}~
\sum_{j=1}^{e}\sum_{k=1}^{r_j}a_{j,k}p^{e-j}\equiv 0~(mod~p^e).
\end{equation}

Suppose the signature of $\Pi$ is $(g~|~-).$ Then the Riemann-Hurwitz
formula (\ref{RH1}) gives
\begin{equation}\displaystyle\label{RH3}
2g-2=p^e\left(-2+\sum_{j=1}^{e}r_j\left(1-\frac{1}{p^j}\right)\right).
\end{equation}
Let $H=\mathbb Z_p$ be the cyclic group of order $p$ generated by some
element $S\in \mathbb Z_{p^e}$ of order $p$.
According to Corollary (\ref{fixedpts4}) the number of fixed points of
$S$ is
$$\displaystyle
t=\sum_{j=1}^{e}r_jp^{e-j}.
$$
Using this value of $t$ in Equation (\ref{RH2}) and then comparing 
Equations (\ref{RH2}) and (\ref{RH3}) we get 
$$\displaystyle
k=1-p^{e-1}+\frac{1}{2}\sum_{j=1}^{e}r_jp^{e-j}\left(p^{j-1}-1\right)
=\frac{1}{2}\sum_{j=2}^{e-1}r_jp^{e-j}\left(p^{j-1}-1\right)
+\frac{1}{2}\left(p^{e-1}-1\right)\left(r_e-2\right).
$$

All terms on the right hand side of this equation are non-negative
since $r_e\ge 2.$ Therefore, if the genus of the action is zero, that
is if $k=0$, we conclude that
$$\displaystyle
r_1\ge 0~\mathrm{is~arbitrary},~r_2=\cdots=r_{e-1}=0~\mathrm{and}~r_e=2.
$$
If $e=1$ we interpret this to mean that the signature of the action
is 
$(0~|~\underbrace{p,\ldots,p}_r)$, where $r\ge 2.$ When $e=2$ the 
interpretation
is that the signature is $(0~|~\underbrace{p,\ldots,p}_r,p^2,p^2),$
where $r$ is arbitrary.

The calculation of the genus is a simple consequence of the
Riemann-Hurwitz formula (\ref{RH1}). This concludes the proof of
Theorem (\ref{cyclicgenusactions}).

We can say more about these actions.  Recall that $T\in\mathbb
Z_{p^e}$ is a generator.  To simplify notation let the elliptic
generators of $\Gamma=\Gamma(0|\underbrace{p,\ldots,p}_{r},p^e,p^e)$
be denoted by
\begin{eqnarray*}\displaystyle
X_1,\ldots,X_r &\mathrm{~of}&~\mathrm{~period~}p,\\ 
Y_1,Y_2&\mathrm{~of}&~\mathrm{~period~}   p^e.
\end{eqnarray*}
Then the epimorphism $\theta:\Gamma\to\mathbb Z_{p^e}$ will satisfy
\begin{eqnarray*}\displaystyle
& &\theta(X_j)=T^{a_jp^{e-1}},~\mathrm{where}~a_j\not\equiv 0~(mod~p),~
1\le j\le r,\\
& &\theta(Y_1)=T^{b_1}~\mathrm{and}~\theta(Y_2)=T^{b_2},~\mathrm{where}~
b_1\not\equiv 0~(mod~p)~\mathrm{and}~b_2\not\equiv 0~(mod~p),\\
& &p^{e-1}\sum_{j=1}^{r}a_j+b_1+b_2\equiv 0~(mod~p^e).
\end{eqnarray*}

Now we prove Theorem (\ref{quatgenusactions}). Let $Q$ denote the
generalized quaternion group $Q(2^n)$ and assume $n\ge 4.$ The case
$n=3$ will be considered at the end of this section.

Some elementary facts about $Q$ are:
\begin{enumerate}
\item The elements of $Q$ are $A^j$ and $BA^j,$ $1\le j\le 2^{n-1}.$
\item The orders of the elements of $Q$ are $1,2,4,\ldots,2^{n-1}.$
\item There is a unique element of order $2$, namely $B^2$, and it 
generates the center of $Q.$
\item The elements of order $4$ are $A^{\alpha 2^{n-3}}$ 
 and $BA^a$, where $\alpha$ is odd
 and $a$ is arbitrary.
\item The elements of order $2^j$, where $3\le j\le n-1$, are 
$\displaystyle A^{\alpha 2^{n-1-j}}$, where $\alpha$ is odd.
\end{enumerate}

Now suppose $M$ is a Riemann surface with an action by $Q$ so that
$M/Q=\mathbb P^{1}(\mathbb C).$ For now we do not assume that the
genus of the action is zero. Then there exists a short exact sequence
$ 1\to\Pi\to\Gamma\stackrel{\theta}{\longrightarrow} Q \to 1, $ as in
(\ref{ses1}), and sig$(\Gamma)$ is given by (\ref{quatsig}).

We need some notation in order to characterize the epimorphism $\theta
:\Gamma\to Q .$ Thus let
$$
X_{j,k},~1\le j\le n-1,~1\le k\le r_j,
$$ 
be  elliptic generators of $\Gamma$ chosen to satisfy the relations
\begin{enumerate}
\item $\displaystyle X_{j,k}^{2^j}=1,~1\le j\le n-1,~1\le k\le r_j.$
\item $\displaystyle X_{1,1}\cdots X_{1,r_1}X_{2,1}\cdots X_{2,r_2}
\cdots X_{n-1,1}\cdots X_{n-1,r_{n-1}}=1.$
\end{enumerate}
Then $\theta$ must preserve the order of the elliptic generators and 
therefore we have
\begin{enumerate}
\item $\displaystyle \theta(X_{1,k})=B^2$, $1\le k\le r_1.$
\item $\displaystyle 
\theta\left(X_{2,k}\right)=
\displaystyle A^{\alpha_{2,k}2^{n-3}}$ 
or $\displaystyle BA^{a_{2,k}}$, where $\displaystyle \alpha_{2,k}$ is odd
and $\displaystyle a_{2,k}$  is arbitrary. 
\item $\displaystyle \theta\left(X_{j,k}\right)=A^{\alpha_{j,k}2^{n-1-j}},~
\mathrm{where}~3\le j\le n-1,~ 1\le k\le r_j~
\mathrm{and~the~}\alpha_{j,k}~\mathrm{are~odd}.$
\end{enumerate}
Let the number of $\displaystyle X_{2,k}$ mapping to
elements of the form $BA^{a}$ be $r.$ Then we must have $r>0$ in order
that $\theta$ be onto. Moreover, in order that $\theta$ be well defined
we must have
\begin{eqnarray}\displaystyle\label{mod4eq}
2r_1+r\equiv 0~(mod~4).
\end{eqnarray}
This follows from counting the number of $B$'s in the product relation
obtained by applying $\theta$ to $\displaystyle X_{1,1}\cdots
X_{n-1,r_{n-1}}=1.$ There is also a restriction coming from the power
of $A$, but we do not need it yet.

>From Equation (\ref{mod4eq}) we see that $r$ must be even, and since
$r>0$, we get $r_2\ge r\ge 2.$

Now $Q$ has a unique subgroup $H$ of order $2$, namely $<B^2>$, and
therefore we need only consider the genus of $M/H$ in order to answer
Problem (\ref{genuscond1}). 
Again according to Corollary (\ref{fixedpts4}), the number of fixed points
of $B^2$ is given by
$$\displaystyle 
t=\sum_{j=1}^{n-1}r_j2^{n-j}.
$$

Assume the signature of $\Pi$ is $(g~|~-).$ Then the Riemann-Hurwitz
formulas (\ref{RH1}) and (\ref{RH2}) give
\begin{eqnarray*}\displaystyle
2g-2&=&2^n\left(-2+\sum_{j=1}^{n-1}r_j\left(1-\frac{1}{2^j}\right)\right)\\
    &=&2\left(2k-2+t\left(1-\frac{1}{2}\right)\right)\\
    &=&4k-4+\sum_{j=1}^{n-1}r_j2^{n-j}.
\end{eqnarray*}
Assuming $k=0$ we easily derive the following equation
\begin{eqnarray}\label{quateq}
r_2\left(2^{n-1}-2^{n-2}\right)+ r_3\left(2^{n-1}-2^{n-3}\right)+\cdots+
+ r_{n-1}\left(2^{n-1}-2\right)  =2^n-2.
\end{eqnarray}
It follows that $r_{n-1}$ must be odd, since the right hand side of 
Equation (\ref{quateq}) is congruent
to $-2$ modulo $4.$ But then we must have
$r_{n-1}=1$, for otherwise the left hand side would be greater than 
$2^n-2$ (this uses $r_2\ge 2$).

Now it follows that the only solutions of Equation (\ref{quateq}) are
$$\displaystyle
r_2=2,r_3=\cdots=r_{n-2}=0~
\mathrm{and}~r_{n-1}=1.
$$
If $n=4$ we interpret this to mean
$\displaystyle
r_2=2~
\mathrm{and}~r_3=1.
$

Thus $r=2,$ and then from Equation (\ref{mod4eq}) we see that $r_1$
must be odd. This means that the only possible signatures are
$(0|\underbrace{2,\ldots,2}_{r_1},4,4,2^{n-1}),$ where $r_1$ is odd.
All of these signatures can be realized by genus-zero actions.

If $n=3$ then the above, suitably interpreted, remains valid. In
particular $r_1$ must be odd and $r_2=3.$

The statement about the genus of the action follows easily from the
Riemann-Hurwitz formula (\ref{RH1}).  This concludes the proof of
Theorem (\ref{quatgenusactions}).

To simplify the notation let the  elliptic generators of $\Gamma$ be denoted 
by
\begin{eqnarray*}\displaystyle
X_1,\ldots,X_r~&\mathrm{of~period}&~2\\
Y_1,Y_2~&\mathrm{of~period}&~4\\
Z~&\mathrm{of~period}&~2^{n-1}
\end{eqnarray*}

The homomorphism $\theta :\Gamma\to Q(2^n)$ must be onto, respect the
relations in $\Gamma$ and have a torsion free kernel.  Therefore, for
$n\ge 4:$
\begin{eqnarray*}\displaystyle
& &\theta(X_1)=\cdots=\theta(X_r)=B^2,
\theta(Y_1)=BA^{a_1},\theta(Y_2)=BA^{a_2},
\theta(Z)=A^{\alpha}\\
& &\mathrm{where}~a_1,a_2~\mathrm{ are ~arbitrary},~\alpha~\mathrm{is~odd},~
r~\mathrm{is~odd},~\mathrm{and}~a_2-a_1+\alpha\equiv 0~(mod~2^{n-1}).
\end{eqnarray*}
There are similar restrictions for $n=3.$

\section{CYCLIC GROUPS}\label{cyclicgrps}

In this section we determine all genus-zero actions for all cyclic
groups.  In particular we prove Theorem (\ref{pqactions}).  Every
cyclic group has genus zero since they can act on $\mathbb P^1(\mathbb
C).$ Moreover the cyclic groups $\mathbb Z_{p^e}$, where $p\ge 2$ is
any prime, admit actions of genus zero on Riemann surfaces of
arbitrarily high genus.  See Theorem (\ref{cyclicgenusactions}).  This
will not be true for other cyclic groups.
 
Now we begin the proof of Theorem (\ref{pqactions}).  Thus let
$G=\mathbb Z_{pq}$, where $p,~q$ are distinct primes, and choose the
presentation
$$G=\left<A,B~|~A^p=1,B^q=1, AB=BA\right>.$$

Suppose $G$ acts on a Riemann surface $M$ so that $M/G=\mathbb
P^1(\mathbb C).$ Then there is a short exact sequence $1\to \Pi\to
\Gamma\stackrel{\theta}{\longrightarrow} G \to 1,$ as in (\ref{ses1}),
where $\Pi$ is a torsion free Fuchsian group of signature $(g~|-)$ and
$$
\mathrm{sig}(\Gamma)
=(0~|~\underbrace{p,\ldots,p}_r,\underbrace{q,\ldots,q}_s,
\underbrace{pq,\ldots,pq}_t).
$$
The Riemann-Hurwitz formula (\ref{RH1}) gives
$$
\displaystyle
2g-2
=pq\left(
-2+r\left(1-\frac{1}{p}\right)+s\left(1-\frac{1}{q}\right)
+t\left(1-\frac{1}{pq}\right)
\right).
$$ 
Let the elliptic generators of $\Gamma$ be 
$$ X_1,\ldots,X_r~\mathrm{of~order}~p;~
Y_1,\ldots,Y_s ~\mathrm{of~order}~q;~
Z_1,\ldots,Z_t~\mathrm{of~order}~pq.$$

Let $H$ be the subgroup of $G$ generated by $A.$ Applying Corollary
(\ref{fixedpts3}) we see that there are $rq+t$ fixed points of $A.$
The signature of $\Gamma^{'}=\theta^{-1}(H)$ will have the form
$(k|\underbrace{p,\ldots,p}_{rq+t})$ and therefore by (\ref{RH2}) we
see that
$$
 \displaystyle 
2g-2=p\left(2k-2+(rq+t)\left(1-\frac{1}{p}\right)\right).
$$
If $k=0$ then a comparison of these formulas for $2g-2$ gives $s+t=2.$
The same argument applied to the subgroup generated by $B$ gives
$r+t=2.$ The only solutions of these equations are
$$(r,s,t)=(0,0,2), (1,1,1),(2,2,0).$$
All $3$ solutions yield actions with genus zero. The genus in each case
is computed by the Riemann-Hurwitz formula. This completes the proof
of Theorem (\ref{pqactions}). 

It is simple enough to describe the epimorphism $\theta$ in each of
these 3 cases.
\begin{enumerate}
\item $\theta(Z_1)=A^aB^b$ and $\theta(Z_2)=A^{-a}B^{-b},$ where
$1\le a\le p-1$ and $1\le b\le q-1.$
\item $\theta(X_1)=A^a,~\theta(Y_1)=B^b$ and $\theta(Z_1)=A^{-a}B^{-b},$
where $1\le a\le p-1$ and $1\le b\le q-1.$
\item $\theta(X_1)=A^a,~\theta(X_2)=A^{-a},~
\theta(Y_1)=B^b$ and $ \theta(Y_1)=B^{-b},$ where 
$1\le a\le p-1$ and $1\le b\le q-1.$
\end{enumerate}

\bco\label{cyclicactions}
Suppose $G$ is a cyclic group $\mathbb Z_n$ with at least $3$ distinct
primes dividing $n.$ Then $G$ does not admit a genus-zero action on a
Riemann surface of positive genus.  \eco \bpf If $G$ did admit such an
action so would $\mathbb Z_{pqr},$ where $p,q,r$ are distinct primes
dividing $n.$ According to Theorem (\ref{pqactions}) we would have
$$\displaystyle
2^\alpha(p-1)(q-1)=2^\beta(p-1)(r-1)=2^\gamma(q-1)(r-1),
$$
where $\alpha,\beta,\gamma$ are either $0$ or $-1.$ This is not possible.   
\epf
The next theorem is proved in a similar fashion, but requires a tedious
case-by-case analysis.  We omit the details.
\bth\label{p^2qactions}
$\mathbb Z_{p^2q},$ where $p,q$ are distinct primes, does not admit a 
genus-zero action on a Riemann surface of positive genus.
\end{theorem}
\section{ZASSENHAUS METACYCLIC GROUPS}\label{zassgrps}

In this section we determine which ZM groups $G_{m,n}(r)$ have genus
zero and then prove Theorem (\ref{ZMactions}). First note that
$G_{m,2}(-1),$ where $m$ is odd, is the dihedral group $D_{2m}$ and
therefore acts on $\mathbb P^1(\mathbb C).$ Thus $G_{m,2}(-1),$ for
$m$ odd, has genus zero. By a routine, but lengthy computation, it is
possible to show that the only ZM groups acting on $\mathbb
P^1(\mathbb C)$ are $G_{m,2}(-1).$
 
Let $G$ be a ZM group $G$ other than $G_{m,2}(-1).$ Then, according to
Theorem (\ref{g=1}), any genus-zero action must be on a surface of
genus $g>1.$ Thus the $pq$-conditions, see Definition (\ref{pqcond}),
must hold.

Therefore a necessary condition for $G$ to admit an action of genus
zero is that $m,n$ have at most $2$ primes in their prime power
factorizations.  See Corollary (\ref{cyclicactions}).  Moreover,
according to Theorem (\ref{p^2qactions}), both $m$ and $n$ must be
either prime powers or a product of distinct primes.

On the other hand, if our goal is to show that a certain $G$ does not
admit an action of genus zero then we may assume $m$ is a prime $p$
since $G$ contains the ZM subgroup generated by $A^k$ and $B$, where
$k$ is chosen so that the order of $A^k$ is $p$. If the subgroup does
not admit a genus-zero action neither will $G.$ There are then $2$
cases to consider.  First we could have $n=q^e$ and secondly we could
have $n=q_1q_2$, where $q_1,q_2$ are distinct primes.

Assume $n=q_1q_2.$ If $d$ is the order of $r$ modulo $p$ then there
are $3$ choices for $d$, namely $d=q_1,~q_2$ or $ q_1q_2.$ The
possibilities for $\displaystyle\frac{n}{d}$ are then $q_2,~q_1$ or
$1$ respectively.  However, the $pq$ conditions are equivalent to the
statement that every prime divisor of $d$ also divides $\displaystyle
\frac{n}{d}.$ Thus there are no choices of $d$ satisfying the $pq$
conditions and so $n=q^e.$

Therefore we first consider the case of a ZM group $G$ where $m=p$ and
$n=q^e$, where $p,q$ are primes. The numerical restrictions on $m,~n$,
see (\ref{ZMgrps}), imply that $p,~q$ are distinct primes.  Moreover
we may assume that $e\ge 2$ since otherwise $G$ would be cyclic. In
fact we will start with the special case of a ZM group $G$ with $m=p$
and $n=q^2$, where $p$ and $q$ are distinct odd primes.  The case
where one of the primes is $2$ will be treated later in this
section. The group $G$ will then satisfy the conditions in Definition
(\ref{pqcond}) if, and only if, $d=q.$

In other words we are considering the group $G$ presented by
\begin{eqnarray}\displaystyle\label{ZMgrp2}
\notag & &\mathrm{generators~:}~A,B\\
& &\mathrm{relations\quad :}~ A^p=1,~B^{q^2}=1,~BAB^{-1}=A^r,~p,q~
\mathrm{distinct~odd~primes}\\ 
\notag & &\mathrm{conditions~:}~GCD((r-1)q,p)=1,~r^q\equiv 1~(mod~p),~\mathrm{and}~
r\not\equiv 1~(mod~p).
\end{eqnarray}
\begin{lemma}
The ZM group in  (\ref{ZMgrp2}) does not admit a genus-zero action.
\end{lemma}
\begin{proof} 
 By 
induction we easily prove that
$$\displaystyle
\left(A^iB^j\right)^k=A^i(1+r^j+r^{2j}+\cdots+r^{(k-1)j})B^{jk}.
$$
Thefore the  possible orders of elements of $G$ are $1,~p,~q,~q^2$ and $pq.$
In particular the order of the element $AB^q$ is $pq$. Let 
$H\cong \mathbb Z_{pq}$ be the subgroup
generated by $AB^q.$

Assume $G$ admits an action of genus zero on some Riemann surface
$M$. Then there exists a short exact sequence $1\to \Pi\to
\Gamma\stackrel{\theta}{\longrightarrow} G\to 1$ as in (\ref{ses1}),
where $\Pi=\Gamma(g~|~-)$ for some $g\ge 0,$ and
$$\displaystyle
\mathrm{sig}(\Gamma)=(0~|~\underbrace{p,\ldots,p}_r,\underbrace{q,\ldots,q}_s,
\underbrace{q^2,\ldots,q^2}_t,\underbrace{pq,\ldots,pq}_u),
$$
for some choice of non-negative integers $r,~s,~t,~u.$ These integers must be 
chosen so that $\theta:\Gamma\to G$ is a well defined epimorphism with 
torsion free kernel.  In particular this means that 
$$
t\ge 2,~\mathrm{and~if}~t=2~\mathrm{then}~r+s+u>0.
$$
The restriction $t\ge 2$ follows from the fact that 
$G_{ab}\cong\mathbb Z_{q^2},$ and therefore there is an epimorphism
$\Gamma\to\mathbb Z_{q^2}.$

The Riemann-Hurwitz formula becomes
\begin{eqnarray*}\displaystyle
2g-2&=&pq^2
\left(
-2+r\left(1-\frac{1}{p}\right)+s\left(1-\frac{1}{q}\right)
+t\left(1-\frac{1}{q^2}\right)+u\left(1-\frac{1}{pq}\right)
\right)\\
&=&-2pq^2+rq^2(p-1)+spq(q-1)+tp(q^2-1)+uq(pq-1)\\
&=&(t-2)pq^2+(r+s+u)pq^2-rq^2-spq-tp-uq.
\end{eqnarray*}

The action of the subgroup $H\cong \mathbb Z_{pq}$ on the Riemann
surface $M$ has genus-zero and therefore, according to Theorem
(\ref{pqactions}), there are only $3$ possibilities for $g$:
$$\displaystyle
g=0,~g=\frac{1}{2}(p-1)(q-1)~\mathrm{or}~ g=(p-1)(q-1).
$$
We will analyze each case separately and conclude that the group 
$G$ does not admit an action of genus zero.

\begin{enumerate}
\item Assume $g=0.$ Then necessarily
$$\displaystyle 
r\left(1-\frac{1}{p}\right)+s\left(1-\frac{1}{q}\right)
+t\left(1-\frac{1}{q^2}\right)+u\left(1-\frac{1}{pq}\right)<2.
$$
The left hand side will be larger than what we get by putting $p=q=3$
and using $t=2.$ This leads to the inequality $6r+6s+8u<2$, and so
$r=s=u=0.$ But then $\theta$ will not be an epimorphism.

\item Assume $\displaystyle g=\frac{1}{2}(p-1)(q-1).$ Then the 
Riemann-Hurwitz equation becomes
$$
(t-2)pq^2+(r+s+u)pq^2-rq^2-spq-tp-uq=(p-1)(q-1)-2,
$$ 
which can be rewritten in the form
$$ (t-2)p(q^2-1)+rq^2(p-1)+(s+1)pq(q-1)+(u-1)q(pq-1)=p-1.$$
>From this equation it follows that $u=0$, for otherwise the left hand side
would be greater than $p-1.$ Recall that $t\ge 2.$ Thus the equation becomes
$$
(t-2)p(q^2-1)+rq^2(p-1)+(s+1)pq(q-1)-q(pq-1)=p-1,
$$
which is equivalent to
$$
(t-2)p(q^2-1)+rq^2(p-1)+spq(q-1)=(p-1)(q+1).
$$
Now it is simple to see that there are no solutions.
\item The last case to consider is $g=(p-1)(q-1).$ By arguments similar
to the second case there are no solutions of the Riemann-Hurwitz equation.
\end{enumerate}
\end{proof}

Next we show that a ZM group of odd order does not admit actions of
genus zero.  Note that by the above considerations we need only prove
this when $m=p$ and $n=q^e$, where $p$ and $q$ are distinct odd primes
and $e>2.$

That is suppose $G$ is the ZM group presented by
\begin{eqnarray}\displaystyle\label{oddorderZM}
\notag & &\mathrm{generators~:}~A,B\\
& &\mathrm{relations\quad :}~ A^p=1,~B^{q^e}=1,~BAB^{-1}=A^r,~e>2\\
\notag & &\mathrm{conditions~:}~GCD((r-1)q,p)=1,r^{q^e}\equiv 1~(mod~p).
\end{eqnarray} 
\begin{lemma}
A ZM group of odd order does not admit genus-zero actions.
\end{lemma}
\begin{proof}
We need only show that the group $G$ presented in (\ref{oddorderZM}) does not
have genus zero.
The order $d$ of $r$ modulo $p$ will be $q^f$ for some $f$, $1\le f<e.$
If $f=1$ then  $A$ commutes with $B^q,\ldots,B^{q^{e-1}},$  and since $e>2,$ 
it follows that $G$ has a subgroup isomorphic to $\mathbb Z_{pq^2},$ namely 
the subgroup $<A,B^q>.$ But this contradicts Theorem (\ref{p^2qactions}).

Thus suppose $f>1.$ Then consider the subgroup $G_1$ generated by $A$
and $B_1=B^{q^{f-1}}.$ This is a ZM group with presentation
$$\displaystyle
A^p=1,~B_{1}^{q^{e-f+1}}=1,~ B_1AB_1^{-1}=A^{r_1},~
\mathrm{where}~r_1=r^{q^{f-1}}.
$$
But now the order of $r_1$ is $q$, and again we arrive at a contradiction.
Therefore a ZM group of odd order does not admit actions of genus
zero.
\end{proof}
The last case to consider is the case of a ZM group $G$ of even order.
Thus $G$ has a presentation as in (\ref{ZMpres}), where one of $m,n$
is even.  In fact the conditions in (\ref{ZMpres}) imply that $m$ is
odd and $n$ is even. According to Corollary (\ref{cyclicactions}) and
Theorem (\ref{p^2qactions}) the only possibilities for $m,~n$ are
\begin{eqnarray*}
& &
m=p^e~\mathrm{or}~pq,~\mathrm{where}~p,~q~\mathrm{are~odd~primes~and}~p\ne q,\\
& &n=2^f~\mathrm{or}~2q'~\mathrm{where}~q'~\mathrm{is~an~odd~prime}. 
\end{eqnarray*}
If $n=2q'$ then $d=2,~q'$ or $2q'.$ In all cases the $pq$ conditions
can not be satisfied. Therefore $n=2^f.$ Moreover we must have $f\ge 2$ 
to satisfy the $pq$ conditions. In fact $d=2^j$ for some $j,$ 
 $1\le j\le f-1.$

Now assume $m=pq.$ Then the subgroup generated by $A$ and
$B^{2^{f-1}}$ is cyclic of order $2pq.$ This contradicts Corollary
(\ref{cyclicactions}) and so $m=p^e.$ But now Theorem
(\ref{p^2qactions}) implies that $e=1.$ Finally we must have $f=2$ and
$d=2.$ In other words the only possible ZM group satisfying
(\ref{genuscond1}) is $G_{p,4}(-1).$

Next we show that $G=G_{p,4}(-1)$ admits genus-zero 
actions. The following statements are easy to prove.
\begin{enumerate}
\item The order of $G$ is $4p$ and the orders of the elements of $G$ are 
$1,~2,~4,~p$ and $2p.$
\item There is a unique subgroup of order $2,$ namely the subgroup 
generated by $B^2.$
\item There is a unique subgroup of order $p,$ namely the subgroup
generated by $A.$
\end{enumerate}
Assume there is a genus-zero action $G\times M\to M,$ where $M$ has
genus $g.$  Then there is a short exact sequence
$$1\to \Pi\to \Gamma\stackrel{\theta}{\longrightarrow}
 G\to 1$$
as in (\ref{ses1}), where $\Pi=\Gamma(g~|~-)$ and
$$\displaystyle
\mathrm{sig}(\Gamma)=(0~|~\underbrace{2,\ldots,2}_r,\underbrace{4,\ldots,4}_s,
\underbrace{p,\ldots,p}_t,\underbrace{2p,\ldots,2p}_u).
$$ 
The Riemann-Hurwitz formula becomes
\begin{eqnarray*}\displaystyle
2g-2&=&
4p\left(-2+r\frac{1}{2}+\frac{3s}{4}+t\left(1-\frac{1}{p}\right)+
u\left(1-\frac{1}{2p}\right)\right)\\
&=&-8p+2pr+3ps+4t(p-1)+2u(2p-1).
\end{eqnarray*}
Applying Corollary (\ref{fixedpts3}) we see that
\begin{enumerate}
\item The number of fixed points of $B^2$ is $2pr+ps+2u.$
\item The number of fixed points of $A$ is $4t+2u.$
\end{enumerate}
Now applying the Riemann-Hurwitz formula (\ref{RH2}) to each of these cases
gives
\begin{eqnarray*}\displaystyle
2g-2&=&2\left(-2+(2pr+ps+u)\frac{1}{2}\right)\\
    &=&p\left(-2+(4t+u)(1-\frac{1}{p})\right)
\end{eqnarray*}
Working with these equations we see that the only solution is 
$$r=0,s=2,t=1,u=0.
$$
This completes the proof of Theorem (\ref{ZMactions}).

To describe the epimorphism $\theta:\Gamma\to G$ associated to the genus-zero
action  
let the elliptic generators of $\Gamma$ be $X,~Y$ and $Z,$ chosen 
so that $$X^4=Y^4=Z^p=XYZ=1.$$
Then all admissible epimorphisms $\theta:\Gamma\to G$ are given by
\begin{eqnarray*}
& &\theta (X)=B^{\epsilon}A^a,~\theta (Y)=B^{-\epsilon}A^b,~\theta (Z)=A^c,\\
& &\mathrm{where}~a,b~\mathrm{are~arbitrary},~1\le c\le p-1,\\
& &\mathrm{and}~-a+b+c\equiv 0~(mod~p).
\end{eqnarray*}

We conclude this section with an example illustrating the action when
$a=1,b=0,c=1$ and $\epsilon =1.$  In this case 
the action of $\Gamma$ on the Poincar\'e disc can be described as
follows. Consider the triangle with vertices $E, F$ and $G$ and angles
$\alpha = \pi/p$, $\beta= \pi/4$ and $\beta$ (respectively) centered
at the origin of the hyperbolic disc. Let $X$ and $Y$ be rotations by an angle
$2\pi/4$ about the vertices $F$ and $G$ respectively, and let $Z$ be
rotation about $E$ with angle $2\pi/p$.  A fundamental domain for
$\Gamma$ is a copy of the shaded triangle $(p,4,4)$ together
with its reflection along the side $EG$. See Figure 1.

\vskip 10pt
\centerline{\psfig{file=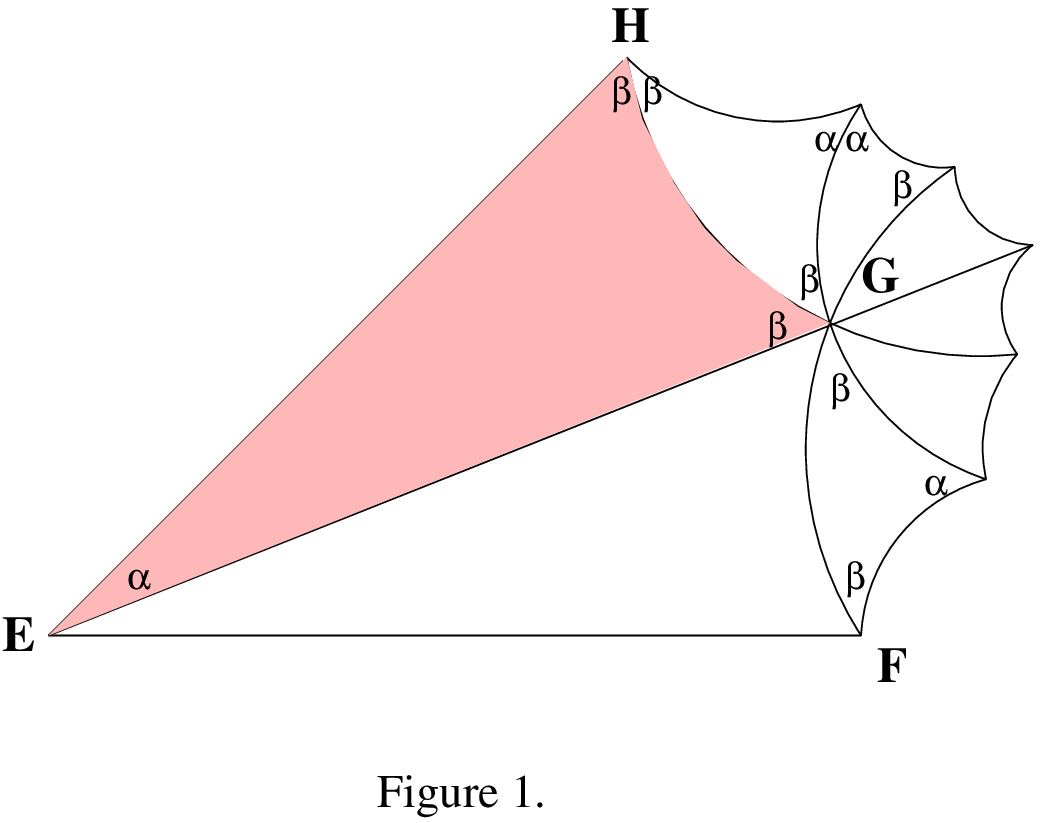, height = 120pt}}
\vskip 10pt

The translates under the rotations $X$, $Y$ and $Z$
of this fundamental domain tessellate the hyperbolic disc. Each rotation
corresponds to the composition of two reflections along a pair of sides. 
Fig. 1  shows four copies of the fundamental domain corresponding
to four rotations about the vertex $G$.

Recall that $G$ acts on a Riemann surface $M$ of genus $p-1$ and one has
the following tower of quotients and group actions
\vskip 10pt
\centerline{\psfig{file=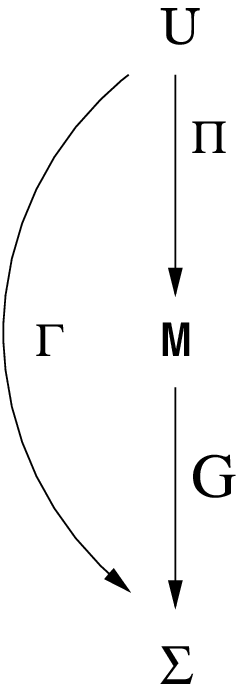, height = 100pt}}
\vskip 10pt
A fundamental domain for $\Pi$ consists of a choice of $4p$ copies of the fundamental
domain of $\Gamma$ (here $|G|=4p$).  An explicit choice can be obtained from the fundamental domain for $\Gamma$ 
(Figure 1) by rotating $p$ times about the vertex $E.$  See Figure 2.
The surface $M$ is obtained from this fundamental domain by making certain
boundary identifications.
\vskip 10pt
\centerline{\psfig{file=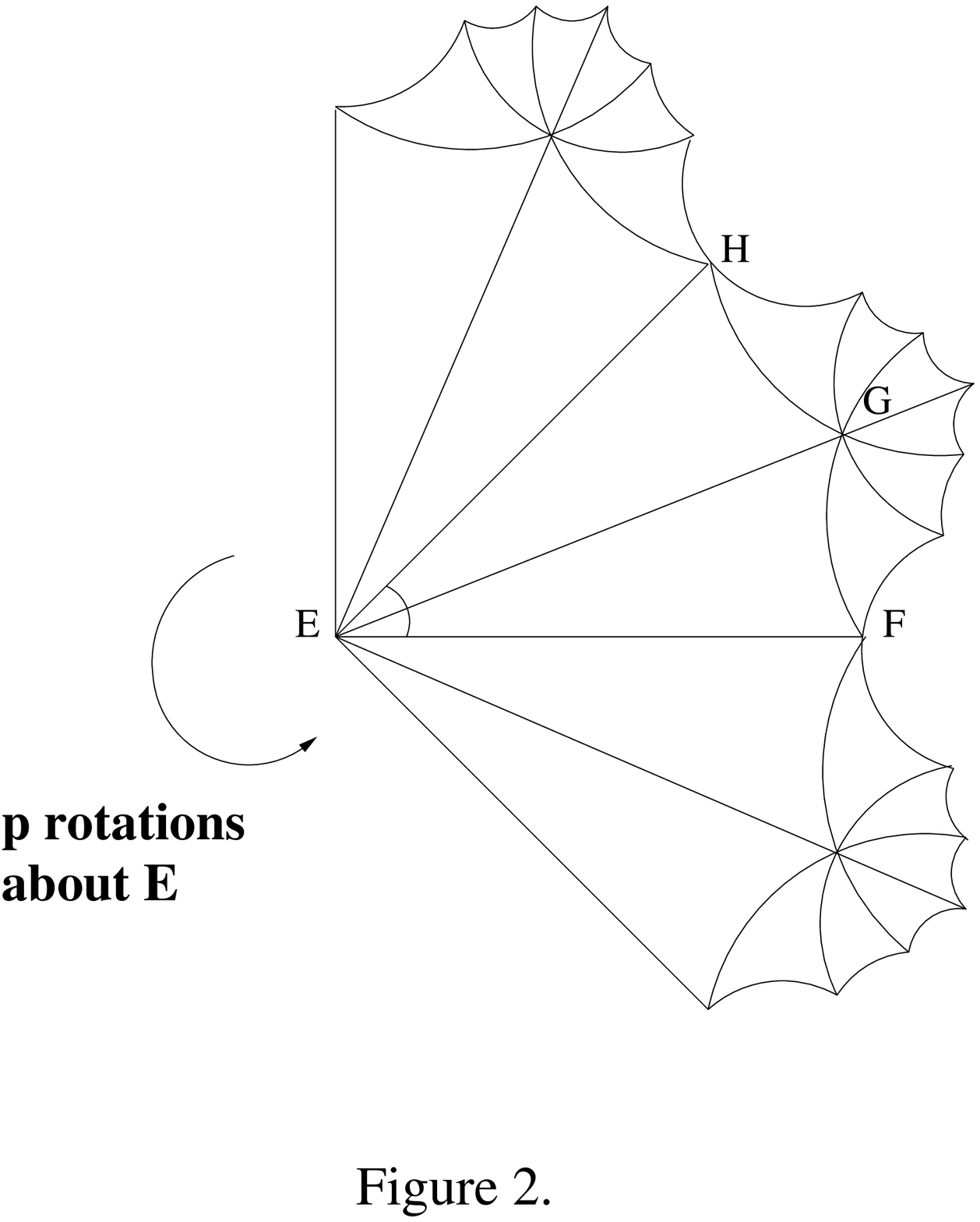, height = 200pt}}
\vskip 10pt

Let $H\subset G$ be the subgroup of order $p$ generated by $A$. Then
there is a short exact sequence 
$1\longrightarrow\Pi\longrightarrow\Gamma'\longrightarrow 
H\longrightarrow 1,$ 
where $[\Gamma :\Gamma'] = [G:H] = 4$ and $\mathbb U/\Gamma' =
M/H$.  A fundamental domain of $\Gamma'$ in $\mathbb U$ will consist
of four copies of the fundamental domain for $\Gamma,$ and this is
depicted in Figure 1. Rotation by $2\pi/p$ about every second vertex in Fig. 1
represents a generator of $H$, see Fig. 3.
The quotient $\mathbb U/\Gamma'$ is obtained
by identifying the sides as depicted, and is easily seen to be $\mathbb P^1(\mathbb C)$.

\vskip 10pt
\centerline{\psfig{file=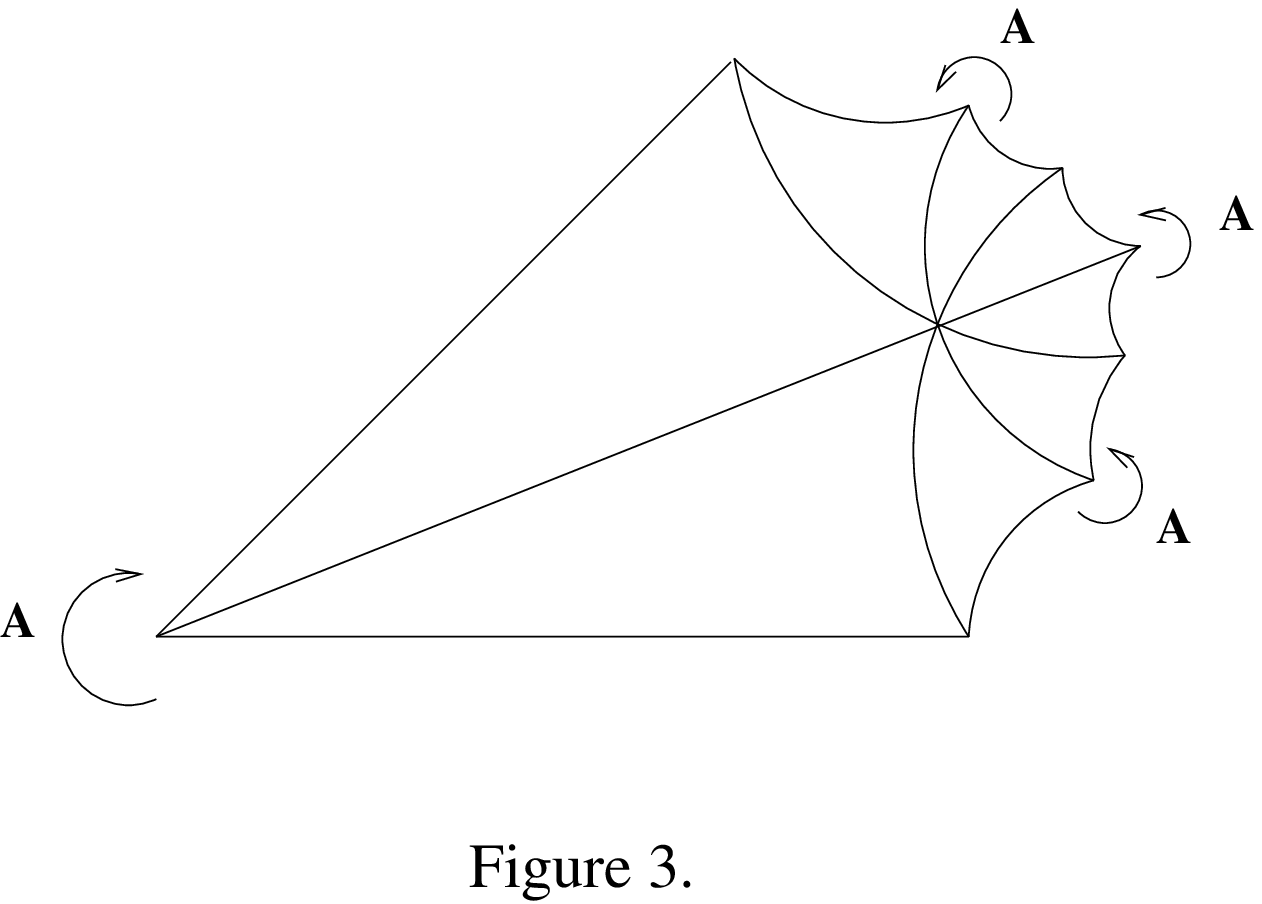, height = 140pt}}
\vskip 10pt

The same analysis can be carried out for the unique subgroup of order
2, $\langle B^2\rangle\subset G$. In this case the fundamental domain
consists of $2p$ copies of the fundamental domain of $\Gamma$ and the
identifications are given as follows (Figure 4).

\vskip 10pt
\centerline{\psfig{file=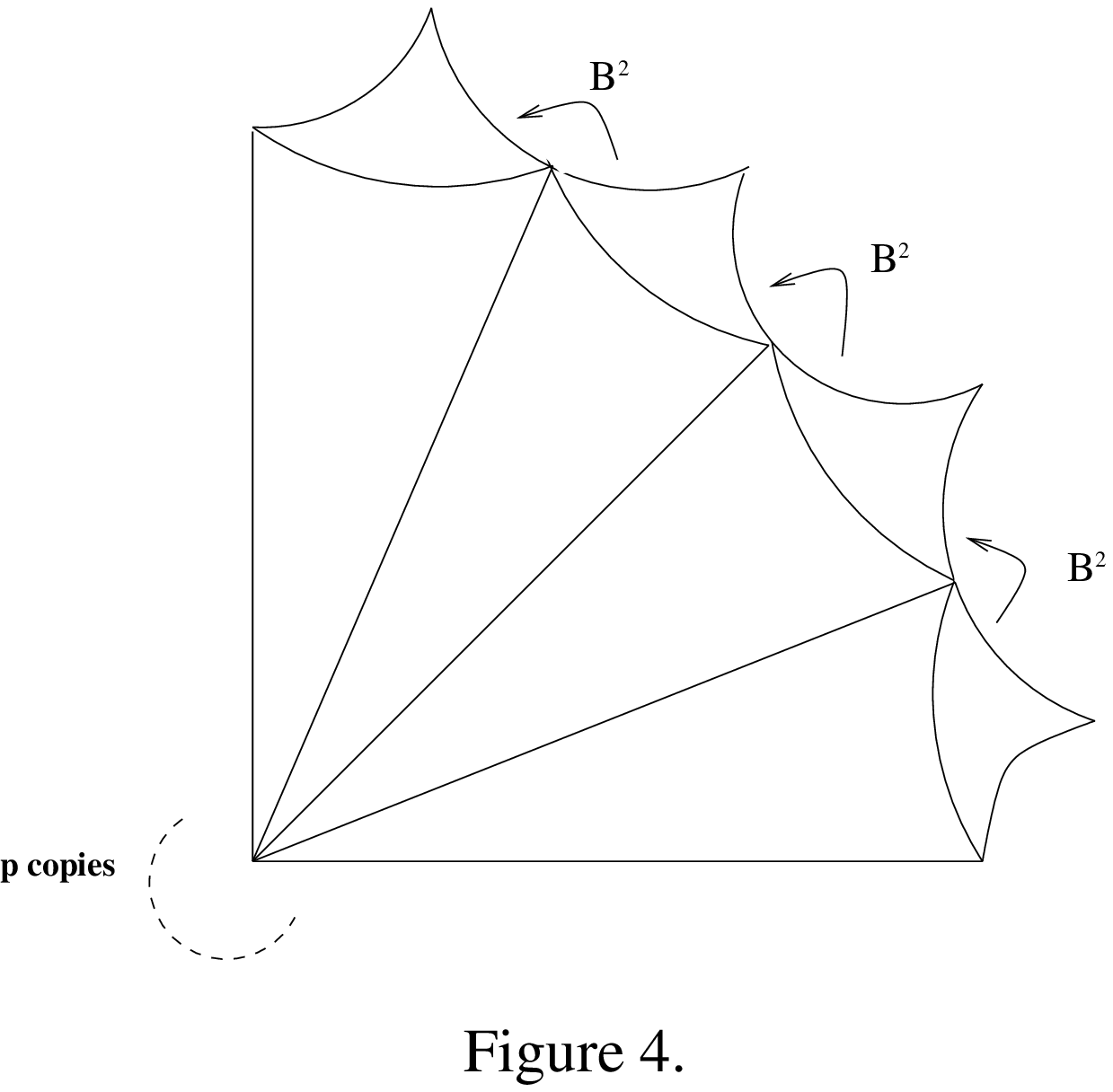, height = 160pt}}
\vskip 10pt

Contiguous pairs of sides on the boundary are identified by $B^2.$
The quotient surface is again $\mathbb P^1(\mathbb C),$
as expected.

\section{THE GENERAL CASE}\label{binaryicos}

In this section we complete our analysis of groups $G$ satisfying
(\ref{genuscond1}) by establishing Theorem (\ref{mainth}). There are
two possibilities to consider, either $G$ is solvable or it is
non-solvable.

First we suppose $G$ is a finite non-solvable group satisfying the
$pq$ conditions, see Definition (\ref{pqcond}).  We want to show that
$G$ does not admit actions of genus zero.  A typical example of such a
group is the binary icosahedral group $I^*$, see Wolf \cite{Wolf}. In
fact any such $G$ must contain $I^*$ as a subgroup and therefore we
need only show that $I^*$ does not admit an action of genus zero.

The binary icosahedral group has order $120$ and contains elements of
orders $1,2,3,4,5,6$ and $10.$ By applying Theorem (\ref{pqactions})
we see that if $I^*$ did admit an action of genus zero on some Riemann
surface $M$ the genus of $M$ would have to be $2.$ Moreover the
signature of the action would be
$$\displaystyle
\mathrm{sig}(\Gamma)=(0~|~\underbrace{2,\ldots,2}_r,\underbrace{3,\ldots,3}_s,
\underbrace{4,\ldots,4}_t,\underbrace{5,\ldots,5}_u,
\underbrace{6,\ldots,6}_v,\underbrace{10,\ldots,10}_w).
$$
Applying the Riemann-Hurwitz Formula (\ref{RH1}) with $g=2$
gives the diophantine equation
$$\displaystyle
30r+40s+45t+48u+50v+54w=121.
$$
It is routine to show that this equation does not have any non-negative
integral solutions. 

This proves the following theorem.
\begin{theorem}
If $G$ admits an action of genus zero then $G$ is a solvable
 group.
\end{theorem}

Now suppose $G$ is a solvable group satisfying the $pq$ conditions.
According to the table on page 179 of \cite{Wolf} there are $4$ types
of such groups $G$, ( denoted \uppercase\expandafter{\romannumeral1},
\uppercase\expandafter{\romannumeral2},
\uppercase\expandafter{\romannumeral3} and
\uppercase\expandafter{\romannumeral4}), and every type contains some
ZM subgroup $G_{m,n}(r).$ According to Theorem (\ref{ZMactions}) $m$
must be an odd prime $p,$ $n=4$ and $r=-1$ if there is to be a
genus-zero action. This rules out types
\uppercase\expandafter{\romannumeral3} and
\uppercase\expandafter{\romannumeral4}.

Type \uppercase\expandafter{\romannumeral1} is just the ZM
group $G_{p,4}(-1)$ and type \uppercase\expandafter{\romannumeral2} is the 
group $G$ with the following presentation:

\begin{eqnarray*}\displaystyle
\mathrm{generators:}& &A,B,R;\\
\mathrm{~relations:}& & 
A^p=1,~B^4=1,~BAB^{-1}=A^{-1},\\
& &R^2=B^2,~RAR^{-1}=A^{l},~RBR^{-1}=B^{-1};\\
\mathrm{conditions:}& &l\equiv \pm 1~(mod~p).
\end{eqnarray*}
In fact we must have $l=-1$ since otherwise the subgroup $\left<~A,R~\right>$ 
would be cyclic of order $4p,$ contradicting Theorem (\ref{p^2qactions}). 
But now $(RAB)^4=A^4$ and therfore $RAB$ has order $4p$.  Again this is a 
contradiction.

This completes the proof of Theorem (\ref{mainth}).

\vskip 30pt
\flushleft{
Sadok Kallel\\
Universit\'e des Sciences et Technologies de Lille\\
Laboratoire AGAT
U.F.R de Math\'ematiques\\
59655 Villeneuve d'Ascq, France\\
{\sc Email:}  sadok.kallel@agat.univ-lille1.fr

Denis Sjerve\\
Department of Mathematics \\
\#121-1884 Mathematics Road\linebreak
University of British Columbia\\
Vancouver V6T 1Z2\\
{\sc Email:}  sjer@math.ubc.ca}

\end{document}